\documentstyle{amsppt}
\input epsf

\magnification=\magstep1
\voffset = 0.5 true in
\vsize=9 true in
\hsize=6.5 true in
\topmatter
\title
Deformations of Maass forms
\endtitle
\author
D.W.~Farmer \\
S.~Lemurell
\endauthor
\thanks Research of the first author supported in
part by the National Science Foundation and
the American Institute of Mathematics .  
Research of the second author supported in 
part by ``Stiftelsen f\"or internationalisering av
h\"ogre utbildning och forskning'' (STINT).
\endthanks
\thanks
\endthanks
\thanks
\endthanks
\address
American Institute of Mathematics
360 Portage Ave
Palo Alto, CA 94307
farmer\@aimath.org
\endaddress
\address
Chalmers University of Technology, SE-412 96
G\"oteborg, Sweden
sj\@math.chalmers.se
\endaddress           
\abstract
We describe numerical calculations which examine the Phillips-Sarnak
conjecture concerning the disappearance of cusp forms on a noncompact
finite volume Riemann surface $S$ under deformation of the surface. 
Our calculations indicate that if the Teichm\"uller space of $S$ is not trivial
then each cusp form has a set of deformations under which either 
the 
cusp form remains a cusp form, or else it dissolves into a resonance
whose constant term is uniformly a factor of $10^{8}$ smaller than
a typical Fourier coefficient of the form.
We give explicit examples of those deformations in several cases.

\endabstract
\endtopmatter
\document

\NoBlackBoxes


\def\R {{\Bbb R}}
\def\Z {{\Bbb Z}}
           
\def\H {{\Cal H}}

\def\ve{\varepsilon}

\def\({\left(}
\def\){\right)}

\def \mm#1#2#3#4 {\pmatrix \mathstrut#1 & \mathstrut#2 \cr 
			\mathstrut#3 & \mathstrut#4\endpmatrix}

\head
1.  Introduction and statement of results
\endhead

We summarize the basic facts about eigenvalues of the 
Laplacian on compact and noncompact surfaces, and then
describe our calculations.

\subhead 1.1 Weyl's law 
\endsubhead

The Laplacian $\Delta$ on a compact Riemann surface $S$ has
a discrete spectrum $0=\lambda_0<\lambda_1\le \lambda_2\cdots$.
As is traditional, we  write $\lambda_j=\frac14+i t_j^2 = \frac14+R_j^2$,
and we occasionally refer to $R_j$ as the ``eigenvalue.''
There is a precise estimate for the magnitude of $\lambda_n$ 
given by Weyl's law:
$$
N(T):=\#\{|t_n|\le T\} \sim \frac{V}{4\pi} T^2,
$$
where $V=Vol(S)$ and eigenvalues are repeated according to their
multiplicity.

\subhead 1.2 Noncompact surfaces
\endsubhead

If $S$ is noncompact, but has finite volume, then the situation is
more subtle,  for $\Delta$ will have both a discrete and a
continuous spectrum.  See \cite{Iw} for a complete discussion.
The continuous spectrum consists of
the interval $[\frac14, \infty)$, and for $s\in [\frac14, \infty)$
we write $s=\frac14+i t^2$.  The eigenfunctions in the continuous
spectrum are given by Eisenstein series, and the eigenfunctions
for the discrete spectrum are called {\it Maass forms}
or {\it nonholomorphic cusp forms}.
Here Weyl's law states
$$
M(T) + N(T) \sim \frac{V}{4\pi} T^2 ,
$$
where $M(T)$ is the contribution of the continuous spectrum,
given by
$$
M(T) = \frac{1}{4\pi}\int\limits_{-T}^T 
-\frac{\varphi'}{\varphi}\left(\frac12 + it\right) dt ,
$$
where $\varphi$ is the determinant of the scattering
matrix for the Eisenstein series.  

It is a fundamental problem to determine which of 
$M(T)$ or $N(T)$ makes the main contribution to the spectrum.
If $S$ is a surface corresponding to a congruence subgroup,
then techniques from analytic number theory can be used to show that 
$M(T)\ll T\log T$, so then
$$
N(T) \sim \frac{V}{4\pi} T^2 ,
$$
as in the compact case.  
This is also conjectured to hold for any 
arithmetic surface \cite{S1,S2} (see \cite{Ta} for a discussion of
arithmetic Fuchsian groups).
For arithmetic surfaces we also have
the deep conjecture of Selberg that
$\lambda_1\ge \frac 14$.

\subhead 1.3 Nonarithmetic surfaces and the destruction of cusp forms
\endsubhead

If $S$ is a noncompact nonarithmetic finite volume surface,
then work of Colin de Verdiere \cite{C1, C2} and 
Phillips and Sarnak \cite{PS1, PS2, PS3}
suggests that $S$ should in general have no discrete spectrum.
Thus, for noncompact surfaces Maass forms are rare and unusual and
are mostly confined to arithmetic surfaces.

The nonexistence of cusp forms on generic noncompact surfaces is
commonly stated in terms of the ``destruction'' of cusp forms.
Namely, if $S$ has a Maass form with eigenvalue $\lambda$, then for almost all
small deformations of ~$S$, the new surface will have no discrete
eigenvalues in a neighborhood of~$\lambda$.
That is, almost all deformations destroy the cusp form.
In the case of variable negative curvature \cite{C1, C2}
the deformation space is infinite dimensional and a generic perturbation
of the metric will destroy the cusp form.
In the constant negative curvature case \cite{PS1, PS2, PS3}
deformations are given by the Teichm\"uller space of the surface,
which is finite dimensional.  This situation is much more subtle.

As described by Phillips and Sarnak, under deformations of the surface
the elements of the discrete spectrum tend to become resonances,
that is, poles of the scattering matrix of the  Eisenstein series.
In the case of a group with one cusp, the Eisenstein series satisfies 
a functional equation of the form $E(z; s)=\varphi(s)E(z; 1-s)$,
where $\varphi(s)\varphi(1-s)=1$.  A cusp form is destroyed if its
eigenvalue moves off the line $\Re(s)=\frac12$ to become a pole
of $\varphi(s)$ in $\Re(s)<\frac12$.  
See also \cite{L,P,W}.

Note that by the expression for $M(T)$ in terms of $\varphi(s)$ above,
this process does not change Weyl's law, although it does change
the balance between $M(T)$ and $N(T)$.

The pole of $\varphi(s)$ in $\Re(s)<\frac12$ corresponds to a zero of
$\varphi(s)$ in $\Re(s)>\frac12$.  H.~Avelin\cite{A} has directly
verified the destruction of cusp forms by tracking such zeros.
In this paper we take a complementary approach and directly track
the cusp forms along deformations which do {\it not\/} destroy
the form.  It would be interesting to compare our calculations
with the results of Avelin.

\subhead 1.4 The results of our calculations
\endsubhead

In this paper we describe numerical calculations 
on noncompact surfaces $S$ 
which indicate that
while it is true that almost all deformations destroy any given cusp form, 
individual cusp forms have deformations along which they are not destroyed.
Our calculations indicate that, 
if $S$ has a nontrivial Teichm\"uller
space, then for each Maass form on $S$ there is a continuous family
of Teichm\"uller deformations on which (deformations of the) Maass form lives.
We explicitly describe these deformations in some simple cases.

A summary of our observations is as follows.  

\item{1.} Suppose $S_0$ has one cusp, the Teichm\"uller space of $S_0$
has dimension~$d$, and $f_0$ is a Maass form on~$S_0$
with eigenvalue~$\lambda_0$.  
Then there is a continuous $d-1$ parameter
family of deformations $S(t)$, Maass forms
$f(t)$, and a function $\lambda(t)$,
such that $S(0)=S_{0}$, $f(0)=f_0$, and $\lambda(0)=\lambda_0$. 

In the above situation we say that ``$f$ lives on a $d-1$ parameter
family of deformations of~$S$.''  Also, if $f$ and $g$ are Maass
forms on surfaces $S_1$ and $S_2$, respectively,
then we say that $f$ and $g$ are {\it equivalent} if 
there exists a continuous family of deformations
which sends $S_1$ to $S_2$ and sends $f$ to $g$.

Note that these calculations are reasonable in terms of
the Phillips-Sarnak phenomenon.  The condition that the
eigenvalue remain on the line $\Re(s)=\frac12$ 
imposes one condition on the deformation, and that condition should be
satisfied on a codimension~1 subset of Teichm\"uller space.
Furthermore, that set should be a real analytic subvariety.

\item{2.} There exist Maass forms $f$ on a surface $S$ with two
dimensional Teichm\"uller space, such that $f$ lives on 
two independent $1$~parameter families of deformations of~$S$.

\item{3.} There exist Maass forms which are not equivalent
to a Maass form on any arithmetic surface.

Since the definition of {\it cusp form} requires invariance under a group 
and vanishing constant term in the Fourier expansion, neither of
those conditions can be verified by a floating point numerical
calculation.  In particular, the functions which we claim are cusp forms
may not actually be invariant under the group, and more critically to 
our purpose here, may not actually vanish at the cusp.  As we discuss
in Section~4 (see particularly Section 4.1.1), our calculations 
are done to a precision of one part in $10^8$, so there is no way for
us to rule out the possibility that our ``cusp forms'' actually have
a nonzero constant term which is a factor of $10^8$ smaller than the
other Fourier coefficients.  Throughout the paper we speak as if
we are finding cusp forms, but it is possible that we are
actually finding one-parameter families of non-cusp forms having
unusually small constant Fourier coefficient.

In the next section we give basic definitions and
describe our calculations.  In the following section
we describe two 
families of Riemann surfaces
for which we explicitly found deformations which do not
destroy a cusp form.  In the final section we 
describe 
our calculations.

\head
2. Definitions and description of calculations
\endhead

We summarize standard facts about the relationship between 
Riemann surfaces and subgroups of $PSL(2,\Bbb R)$, 
and we describe the surfaces on which we performed our calculations.

\subhead 2.1 Subgroups of $PSL(2,\Bbb R)$ and $\H/\Gamma$ 
\endsubhead

The group $PSL(2,\Bbb R)$ acts on the upper half--plane
$\H=\{z=x+iy\ :\ x,y\in \Bbb R,\, y>0\}$ by linear fractional
transformations:
$$
\mm abcd (z)=\frac{az+b}{cz+d}.
$$
Here $\H$ is equipped with the hyperbolic metric
$ds^2=y^{-2}(dx^2+dy^2)$ and area element $dA=y^{-2}dx\,dy$.
If $\Gamma\subset  PSL(2,\Bbb R)$ is a discrete subgroup then
$\H/\Gamma$ is a Riemann surface of constant negative curvature~$-1$.  
Throughout the paper we consider the case where
$\H/\Gamma$ is noncompact and has finite area, and in this case
we say that $\Gamma$ is a cofinite subgroup of~$PSL(2,\Bbb R)$.
The surface $\H/\Gamma$ can be compactified by the addition of
a finite set of points.  These missing points are referred to
as ``cusps'', and they correspond to (conjugacy classes of)
parabolic subgroups of~$\Gamma$.

We use the standard notation for the signature of  Riemann surface
$S$:
$$
sig(S)=\{g,\{m_1,...,m_k\},\nu\},
$$
where $g$ is the genus of $S$, $m_1$,...,$m_k$ are the orders of
the elliptic points of $S$, and $\nu$ is the number of cusps of~$S$.
In the above notation, the Teichm\"uller space of $S$ has (real)
dimension $6g-6+2k+2\nu$.

If $sig(S)$ is as given above then $S=\H/G$ where $G\subset PSL(2,\R)$ has presentation
$$
G=\langle a_1,b_1,...,a_g,b_g,
e_1,...,e_k,p_1,...,p_\nu\ | \ 
	e_j^{m_j}=1, [a_1,b_1]...[a_k,b_k]e_1 ... e_k p_1...p_\nu=1 \rangle .
$$
The equivalence between Riemann surfaces $S$ and groups 
$\Gamma\subset PSL(2,\R)$ allows us to phrase everything in terms
of the group.  This is useful for computer calculations
because we can represent the group by a convenient set of generators.

\subhead 2.2 Maass forms and Hecke congruence subgroups
\endsubhead

A {\it Maass form} on 
a group $\Gamma\subset PSL(2,\R)$ is a function $f:\H\to\R$ which satisfies:
\item{}(1.1) $f(\gamma z)=f(z)$ for all $\gamma\in\Gamma$,
\item{}(1.2) $f$ vanishes at the cusps of $\Gamma$, and
\item{}(1.3) $\Delta f = \lambda f$ for some~$\lambda>0$,

\noindent where 
$$
\Delta=-y^2\(\frac{\partial^2}{\partial x^2} 
	+ \frac{\partial^2}{\partial y^2}\) 
$$
is the Laplace--Beltrami operator on~$\H$.
Note that Maass forms on $\Gamma$ are elements of the discrete spectrum of
$\Delta$ on $S=\H/\Gamma$.

If $\Gamma$ is cofinite then it contains a parabolic element,
which we can conjugate to $T=\mm 1101 $.  Thus, we may assume
$f(z)=f(z+1)$.  By (1.3) we find that $f(z)$ has a Fourier 
expansion of the form
$$
f(z)=\sqrt{y}\sum_{n\not=0} a_n  K_{iR}(2 \pi |n| y)\exp(2 \pi i n x),
$$
where $K_\nu(t)$ is a Bessel function 
and $\lambda=\frac14+R^2$.

Maass forms naturally arise in number theory 
\cite{Iw} in the case of
$\Gamma=\Gamma_0(q)$.  Here $\Gamma_0(q)$ is the Hecke congruence
group
$$
\Gamma_0(q)=\left\{\mm abcd \in PSL(2,\Z) \ :\ q|c \right\}.
$$
It is usually more convenient to consider a slightly 
larger group $\Gamma^*_0(q)$, 
defined as follows.  Suppose $t|q$ and $(t,q/t)=1$, and choose $a,b$
so that $at-bn/t=1$.
Let
$$
H_t(q)=\mm {at}{b}{q}{t} 
$$
and note that $H_t(q)$ normalizes $\Gamma_0(q)$, the $H_t(q)$ are
defined up to multiplication by an element of $\Gamma_0(q)$,
and 
$H_t(q)^2 \in \Gamma_0(q)$.  The $H_t(q)$ are called the 
{\it Fricke involutions}
for $\Gamma_0(q)$.  If $q$ is squarefree then we let 
$\Gamma^*_0(q)$ be the group generated by $\Gamma_0(q)$ and all
of the $H_t(q)$ for $t|q$.  This is a maximal discrete subgroup
of $PSL(2,\R)$.  If $q$ is not squarefree then 
the group generated by $\Gamma_0(q)$ and all
of the $H_t(q)$ for $t||q$ may or may not be a maximal discrete group.
If it is then we call it $\Gamma^*_0(q)$.  If it isn't then we adjoin
as many more involutions as possible.  See below for examples with
$q=8$ and~$9$.  We write $S_0(q)=\H/\Gamma_0(q)$ and
$S_0^*(q)=\H/\Gamma_0^*(q)$.
See \cite{Co} for an interesting discussion 
of the relationship between $S_0(q)$ and $S_0^*(q)$.

The groups $\Gamma^*_0(q)$ are more convenient for our
purposes because $\Gamma^*_0(q)$ has
fewer cusps than $\Gamma_0(q)$, and this simplifies
the search for cusp forms.  In fact, all of the cases we consider here
have just one cusp.  Since the Teichm\"uller
spaces of $\Gamma_0(q)$ and $\Gamma^*_0(q)$ are essentially the same, there
is no loss in considering the larger group.

If $q=1,2,3$, or $4$ then $\Gamma_0^*(q)$ is conjugate to a
Hecke triangle group.  These groups cannot be deformed.
Maass forms on Hecke triangle groups have
been extensively studied by Hejhal \cite{H1, H2, H3}.
If $q\ge 5$ then $\Gamma_0(q)$ can be deformed; we concentrate on the cases 
$q=5$, $6$, $8$, $9$, and $11$.

Table 1 gives the signature of $S_0(q)$ and $S_0^*(q)$ for small $q$.

\vskip .2in

\hskip 1in{\vbox{\offinterlineskip
\hrule
\halign{&\vrule#&\strut\quad#\hfil\quad\cr
height2pt&\omit&&\omit&&\omit&\cr
&\hfil $q$ \hfil&
&\hfil $sig(S_0(q))$ \hfil&
&\hfil $sig(S_0^*(q))$ \hfil&
\cr
\noalign{\hrule}
& \hfil 5$\mathstrut$  \hfil &
& \{0,\{2,2\},2\}  \hfil &
& \{0,\{2,2,2\},1\} \hfil &\cr
& \hfil 6  \hfil &
& \{0,\{\},4\}  \hfil &
& \{0,\{2,2,2\},1\} \hfil &\cr
& \hfil 7  \hfil &
& \{0,\{3,3\},2\}  \hfil &
& \{0,\{2,2,3\},1\} \hfil &\cr
& \hfil 8  \hfil &
& \{0,\{\},4\}  \hfil &
& \{0,\{2,2,2\},1\} \hfil &\cr
& \hfil 9  \hfil &
& \{0,\{\},4\}  \hfil &
& \{0,\{2,2,2\},1\} \hfil &\cr
& \hfil 10  \hfil &
& \{0,\{2,2\},4\}  \hfil &
&  \{0,\{2,2,4\},1\}  \hfil &\cr
& \hfil 11  \hfil &
& \{1,\{\},2\}  \hfil &
&  \{0,\{2,2,2,2\},1\} \hfil &\cr
height2pt&\omit&&\omit&&\omit&\cr}
\hrule}\hfil}

\vskip .2in

Note that $S_0^*(5)$, $S_0^*(6)$, $S_0^*(8)$, and  $S_0^*(9)$
have the same signature, so these surfaces are deformations of each other.
That family of surfaces, along with deformations of
$S_0^*(11)$, serve as our examples in this paper.
In the next section we give explicit deformations of these surfaces.

\head 
3. The deformations
\endhead

We explicitly describe deformations of the surfaces 
$S_0^*(q)$ for
$q=5$, $6$, $8$, $9$, and $11$.

\subhead 3.1 Generators
\endsubhead

Since these surfaces have genus~0 and one cusp,
they can be realized from groups generated
by elliptic matrices along with the
single parabolic matrix $T=\mm 1101 $.  Let $r_k(x,y)$ be
the matrix which acts as a rotation by $2\pi/k$ around
the point
$x+iy\in \H$:
$$
r_k(x,y)= \mm {c y - s x}{s (y^2 + x^2)}{-s}{c y + s x} ,
$$
where $c=\cos(\pi/k)$ and $s=\sin(\pi/k)$.

We will give an explicit description of our groups
in terms of $T$ and $r_k(x,y)$.  First we describe various groups 
containing free parameters, and then indicate how they specialize to
various $\Gamma_0^*(q)$.

\subhead 3.2 Signature $\{0,\{2,2,2\},1\}$
\endsubhead

Let $\Gamma_{2,2,2}(a,b)=\langle T,g_1,g_2,g_3 \rangle $
with
$$
g_1=r_2(b,1/\sqrt{a})
\ \ \ \ \ \ 
g_2= r_2(x,\sqrt{y})
\ \ \ \ \ \ 
g_3=Tg_1 g_2 ,
$$
where 
$$\eqalign{
x=& 
\frac12 \left({\frac{2}{a} + b + 
Y
}\right)
\cr
y=& 
\frac12 \left({-\frac{4}{a^2} + 
    \left( 1 - \frac{2}{a} + b \right) \,
     \left( -b + 
Y
\right) }\right)
\cr
Y=&\sqrt{\frac{4}{a^2} + b^2}  \ .
\cr
}
$$
One can check that $\Gamma_{2,2,2}$ is well-defined for 
$(a,b)$ in a neighborhood of $\{(a,0)\ |\ a > 4+\ve \}$,
and that $\Gamma_{2,2,2}(a,b)$
satisfies the relations $g_1^2=g_1^2=g_3^2=1$
and $g_1g_2g_3T=1$.  Thus, $\H/\Gamma_{2,2,2}(a,b)$
has signature $\{0,\{2,2,2\},1\}$.  

\vskip 0.2 in
{

\noindent \epsffile{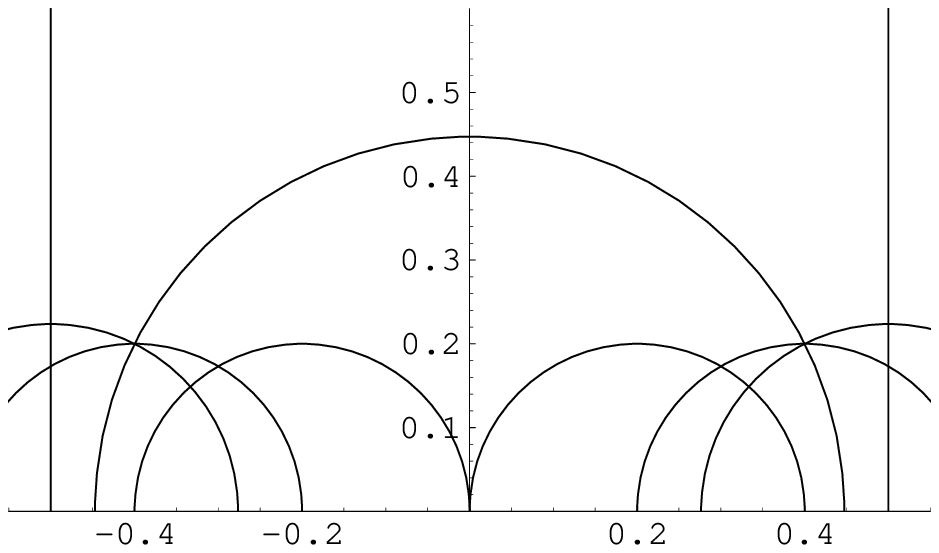}
\hskip 0.3in \epsffile{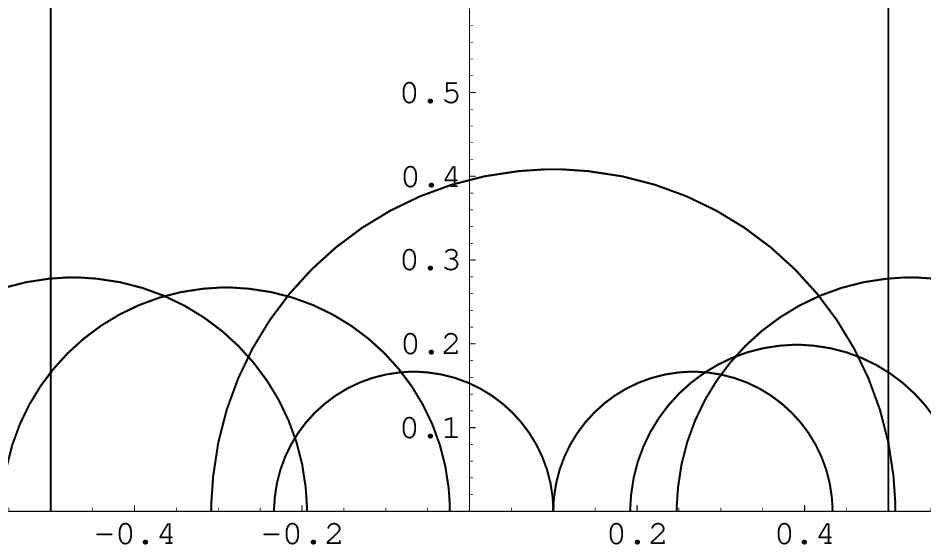}
}

{\sl Figure 3.2.1:  The region above the arcs and between the lines is a fundamental
domain for $\Gamma_{2,2,2}(a,b)$  for $(a,b)=(5,0)$ and 
$(a,b)=(6, \frac1{10} )$, respectively.
} 

\vskip 0.1in

Figure 3.2.1 shows the action of 
the generators of
$\Gamma_{2,2,2}(a,b)$ on the upper half-plane,
along with some additional elements of the group.
Note that the figures have left/right symmetry when $b=0$.

In the figure on the left, the region above the four smallest circles is
a fundamental domain for~$\Gamma_0(5)$, and the region above the three largest
circles is a fundamental domain for~$\Gamma_0^*(5)$.

\subhead 3.3 Signature $\{0,\{2,2,2,2\},1\}$
\endsubhead

Let
$\Gamma_{2,2,2,2}(a,b,c,d)=\langle T, g_1,g_2,g_3,g_4 \rangle $, with
$$
g_1=r_2(a,x)
\ \ \ \ \ \ \ \ \ \ \ \
g_2=r_2(b,y)
\ \ \ \ \ \ \ \ \ \ \ \
g_3=r_2(c,z)
\ \ \ \ \ \ \ \ \ \ \ \
g_4=r_2(1/2,d) ,
$$
where
$$\eqalign{
x=&{\sqrt{\frac{\left(a - b \right)\,
          \left( \left( \frac{1}{2} + a \right)\,
              \left(- \frac{1}{2}  + c \right) + 
                d^2 \right)}{\frac{1}{2} - c}}} \cr
&\cr
y=& d\, {\sqrt{ \frac{\left( a - b \right) \,
          \left( b - c \right) }{\left( \frac{1}{2} + a
            \right) \,\left(\frac{1}{2} - 
            c \right) }  }} \cr
&\cr
z=&{\sqrt{\frac{\left( b - c \right) \,
       \left( \left( \frac{1}{2} + a \right) \,
          \left( - \frac{1}{2}   + c \right)  + 
         d^2 \right) }{{\frac{1}{2} + a}}}}  \ . \cr
}
$$
One can check that $\Gamma_{2,2,2,2}$ 
is well-defined in a neighborhood 
$(-1/3,0,1/3,1/(2\sqrt{11}))$ and that 
$g_1 g_2 g_3 g_4 T=1$.  Thus
$\H/\Gamma_{2,2,2,2}(a,b,c,d)$
has signature $\{0,\{2,2,2,2\},1\}$.

\subhead 3.4 Comparison with $\Gamma_0^*(q)$
\endsubhead

One can check that
$\Gamma_{2,2,2}(v,0)=\Gamma_0^*(v)$
for $v=5$, $6$, or $8$, and
$$
\Gamma_0^*(9)=\mm 1{\frac16}01
        \Gamma_{2,2,2}(9,{\textstyle{\frac 16}})
\mm 1{\frac16}01 ^{-1} .
$$
So $\Gamma_{2,2,2}$ gives the desired deformation of these groups.

Note that there are various symmetries
to $\Gamma_{2,2,2}$ and the
same group can appear several times.  For example,
$\Gamma_0^*(5)$ appears as
$\Gamma_{2,2,2}(5,0)$, $\Gamma_{2,2,2}(20,0)$ and
$\Gamma_{2,2,2}(5,\frac 15)$.  More generally,
$\Gamma_{2,2,2}(a,0)$ is isomorphic to $\Gamma_{2,2,2}(4a/(4-a),0)$.
Also, $\Gamma_{2,2,2}(a,b)=\Gamma_{2,2,2}(a,-b)$.
Other symmetries will be visible in our plots in 
Section 4.

Also, $\Gamma_{2,2,2,2}(-1/3,0,1/3,1/(2\sqrt{11}))=\Gamma_0^*(11)$.
So $\Gamma_{2,2,2,2}$ will be our deformation of $\Gamma_0^*(11)$.

In the next section we use the generators given above 
in our numerical calculations.

\head
4. Tracking the Maass forms
\endhead

Let $\Gamma$ be one of the groups described in the previous section.
We use the methods described in \cite{FJ} to locate 
Maass forms on those groups.  Those methods are summarized in the
next section.  Once one Maass form is located, we search
for a small deformation of the group such that there is a Maass form
with an eigenvalue close to that of the original Maass form.
We then use those two values to interpolate or extrapolate to find 
a starting point for locating more deformations of the original form.  
The end result is a sequence
of closely spaced deformations for which we have Maass forms with slowly
changing eigenvalue and Fourier coefficients.  The great accuracy
to which we find the individual forms, coupled with the 
small changes to the Maass form as we vary the group, provide 
persuasive evidence
that we are indeed tracking the deformation of a Maass form.

\subhead  4.1 Locating Maass forms
\endsubhead

The following is a summary of the methods used to locate an individual
Maass form.  
Given generators $\{g_j\}$ of $\Gamma$, we produce an overdetermined
system of linear equations which uses a truncation of the
Fourier expansion of $f$
$$
\tilde{f}(z)=\sqrt{y}\sum_{|n|\leq M,n\not=0} a_n  K_{iR}(2 \pi |n| y)\exp(2 \pi i n x),
$$
where the $a_i$ are complex unknowns. Note that we assume $a_0=0$
(which excludes all Eisenstein series when we have only one cusp) and
also normalize one of the coefficients (usually $a_1$) to equal~$1$. 
Also note that this of course depends on $R$ (or equivalently on
$\lambda$).  For most of the examples in this paper, we choose~$M$
so that the error caused by the truncation is around $10^{-8}$
for points in the fundamental domain of~$\Gamma$.

We treat the $\{a_n\}$ as $4M-2$
real unknowns. Next we choose $N$
points $z_i$ (where $N>4M$)  on a horizontal line in $\H$. These
points are mapped by the generators to points $g_j z_i=z_i^*$ higher up in
$\H$. If $f$ is a Maass form on $\H/\Gamma$ then $f(z_i)=f(z_i^*)$ (or
more generally $f(z_i)=\chi(g_j)f(z_i^*)$ where $\chi$ is a
character). The $N$ equations 
$$
\tilde{f}(z_i)=\tilde{f}(z_i^*)
$$
constitute an overdetermined system $Ax=b$ in $4M-2$ unknowns. 
If $R$ is an eigenvalue of a Maass form on $\Gamma$, then this 
system should be consistent to within the error caused by the truncation.

Next we determine
the least square solution $\tilde x$ (using $QR$-factorization) to
this system of equations. We then use the norm of the error,
$||A\tilde{x}-b||_2$, as a measure of how close $\lambda$ is to an
eigenvalue. If $\lambda$ is really an eigenvalue, then
$||A\tilde{x}-b||_2$ should be roughly the size of the truncation
error. In our initial calculations in cases where earlier data were
available we found this to be true. We also found that away from
eigenvalues (i.e. if we choose $R$ randomly) the error 
is generally of size~$1$, independent of the size of the truncation error. 
We take $||A\tilde{x}-b||_2$ to be a measure of distance
between $R$ and a ``true'' eigenvalue for~$\Gamma$, and we have
found this measure to vary smoothly and to be consistent 
with various other checks, which we describe below.  Thus, it seems reasonable to 
say that these functions are a factor of $10^8$ closer to 
being invariant under~$\Gamma$ then a randomly chosen function.

There are a number of error checks. If $\Gamma$ is arithmetic
then the Fourier coefficients will be multiplicative, and we
find that the above method produces functions whose coefficients 
are multiplicative to better than  the truncation error.
These can be viewed as independent 1~in~$10^8$ error checks,
which render it very likely that the functions produced are
indeed Maass forms.  In other words, the possibility of a ``false alarm''
is extremely small, and we have high confidence that the program 
is finding the Maass forms for~$\Gamma$.  

For nonarithmetic groups there are no Hecke relations, but there
are other persuasive checks.   We start with a general Fourier 
expansion with complex coefficients.  For the Maass forms we find,
the functions are real to very high accuracy (to an even higher
accuracy then the truncation error).  In general, when we are far from 
an eigenvalue, the system of equations is far from consistent
and the approximate solutions are far from real.

A final check is the size of the Fourier coefficients.  For
arithmetic~$\Gamma$ all the coefficients we have found fit the
Ramanujan-Peterson conjecture~$|a_p|<2$.  For nonarithmetic
groups that bound is not true in general, but it is
still conjectured that~$a_n\ll n^\epsilon$.  And we do in 
fact find that if $R$ is close to an eigenvalue then
the $a_n$ from the least-squares solution are much smaller
(and not growing as a function of $n$)
than those from random~$R$.

Having identified one Maass form, it can serve as
a starting point for locating Maass forms on nearby groups.
If no Maass forms are detected on nearby groups then we can 
state with confidence
that the Maass form is destroyed by all deformations.
If a Maass form is detected, then by finding a succession of 
nearby Maass forms on nearby groups, we can state with confidence
that the Maass form survives as the group is deformed.
We do indeed find continuous families of Maass forms in all cases
which we checked, so we suggest that on all groups which admit
deformations, each Maass form has a continuous family of deformations
on which it lives. 

The following sections illustrate some of the deformations associated
with Maass forms on $\Gamma_{2,2,2}$.  Since the deformation space of this
group is two (real) dimensional, it is possible to visualize the deformations
of the Maass forms.  In the case of $\Gamma_{2,2,2,2}$ there is a 4-dimensional
deformation space, and we must be less direct in our demonstration that
the Maass forms live on a 3-dimensional subspace.

Samples of our data can be found at
\hbox{http://www.math.chalmers.se/$\sim$sj/Maass/}

\subhead
4.1.1 Are we really finding cusp forms?
\endsubhead

One important issue is justifying that we are actually finding cusp forms,
as opposed to the residue of a pole of an Eisenstein series.
By construction our functions vanish at
the cusps (because $a_0=0$ and $\Gamma$ has only one cusp).
But in the nonarithmetic cases it is possible that we are actually
finding a resonance
with extremely small constant term: smaller than the
truncation error.  This would be a
serious concern if we only considered a few examples, and those
involved very small deformations of an arithmetic group.
But, as can be seen below, we consider deformations which
are quite far from an arithmetic group, and we consistently
find our error check to be a factor of $10^8$ smaller than
expected for random data.  So, unless the Maass forms deform
in a way that
whose constant term always stays extremely small
(a factor of $10^8$ smaller than the other coefficients),
then the functions we find are cusp forms.

We have also done calculations where we include a (non-zero) constant
term
$$
\sqrt{y}(a_0\cdot\cos R\ln y+b_0\cdot\sin R\ln y).
$$
We know that there is a one-dimensional space of Eisenstein-series,
since we have exactly one cusp. Numerically we find  a unique solution
at a generic point
in our deformation space. This solution does not
vary (more than the truncation error) when we vary the chosen points
$z_i$. This solution should be the Eisenstein series.  At a point
where we claim that we find cups forms these calculations give
different solutions for each set of points $z_i$. Taking a suitable
linear combination of any two of these solutions gives a solution with
constant coefficients $a_0$ and $b_0$ less than the truncation
error. This solution
(of course) agrees with the one we get when forcing the constant
coefficients to be zero.
Setting $a_0=10^{-5}$ and $b_0=0$ gave errors consistently about $100$ times
larger than when they were both set to zero. We also remark that we
normalized the size of the terms in such a way that the coefficients in
the matrix $A$ were of size $1$ both for the constant
coefficients and the first non-constant coefficients.

It would be interesting to input our deformations into the method
of Avelin~\cite{A} to track the motion of the zeros and poles of the
scattering matrix.  See comment (1) in Section 4.6 for a particularly interesting
example.

If the cusp forms are becoming resonances with uniformly very small
constant term, then it would be interesting to find an explanation
for this phenomenon.

\subhead 4.2 Tracking the eigenvalues
\endsubhead

Our first example is an exhaustive search for all Maass forms 
for $\Gamma_{2,2,2}(a,b)$
in the box
$$
(5\le a\le 6)
\times 
(0\le b\le 0.16) 
\times
(11\le R\le 12)
.
$$
Note that $(a,b)=(5,0)$ corresponds to $\Gamma_0^*(5)$ and 
$(a,b)=(6,0)$ corresponds to $\Gamma_0^*(6)$.

When $b=0$ the group has an extra symmetry, as can be seen in Figure 3.2.1.
In this case the Maass forms are
classified as ``even'' or ``odd'' according to whether they are even
or odd functions of~$x$. 

The plots in Figure 4.2.1 depict the 13 equivalence classes of Maass forms
which intersect the region of the search.
The cube in the upper left shows all of the forms.  The 
left/right, 
front/back,
and up/down axes correspond to $a$, $b$, and $R$, respectively.
The back-left edge corresponds to $\Gamma_0^*(5)$, and the back-right edge
corresponds to $\Gamma_0^*(6)$.  
The back face, where $b=0$, are those deformations
having the extra symmetry described above.  

Odd Maass forms on $\Gamma_0(N)$ survive under 
deformations preserving the left-right symmetry of the fundamental domain
because the Eisenstein series are even and so it is not possible
for the odd forms to dissolve into the continuous spectrum.
This is illustrated in the diagram by the two equivalence classes contained
in the back face of the cube.  Thus, the figure shows that there are
four Maass forms on $\Gamma_0^*(5)$ with eigenvalue
$11\le R\le 12$, with two of them odd and two of them even.

Note:  A separate search found that there were indeed exactly four Maass
forms in that range for $\Gamma_0^*(5)$.  Also, there are even Maass forms
for $\Gamma_0^*(6)$ with $11\le R\le 12$, but for small deformations they
leave the search region.

Proceeding counterclockwise from the cube, the other figures show the projections
onto the back, bottom, and left faces of the cube, respectively.

The view from the top shows an extra symmetry along the curve
from $(a,b)=(5, 0.0854)$ to $(a,b)=(6, 0.1266)$.  We did not trace all of
the Maass forms past their point of symmetry.

Further observations from these figures are made in Section 4.5

\noindent \epsffile{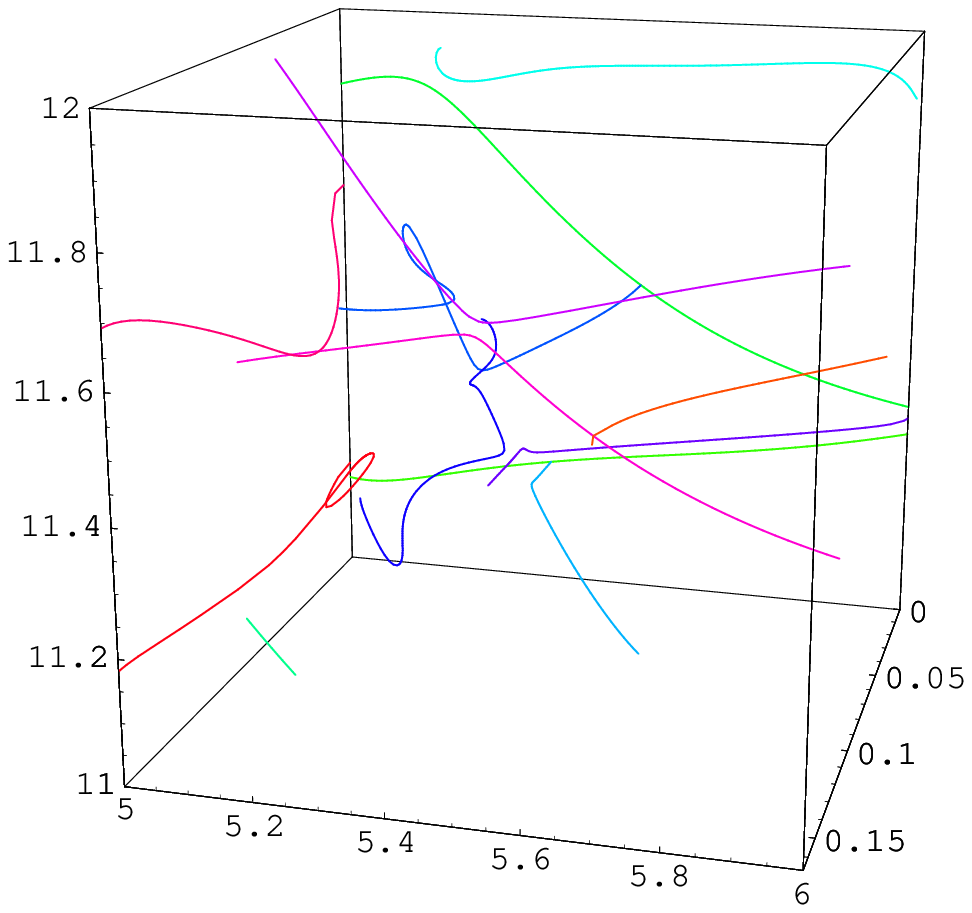}
\hskip .2in \epsffile{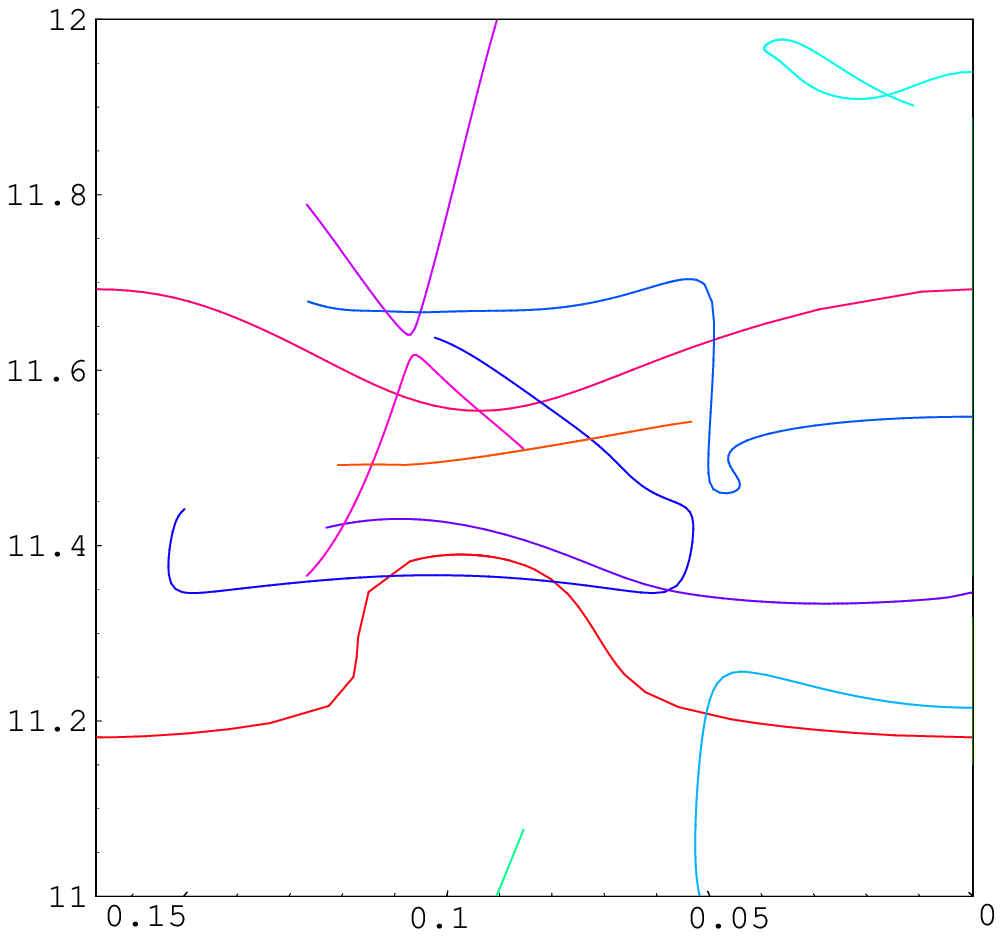}

\noindent \epsffile{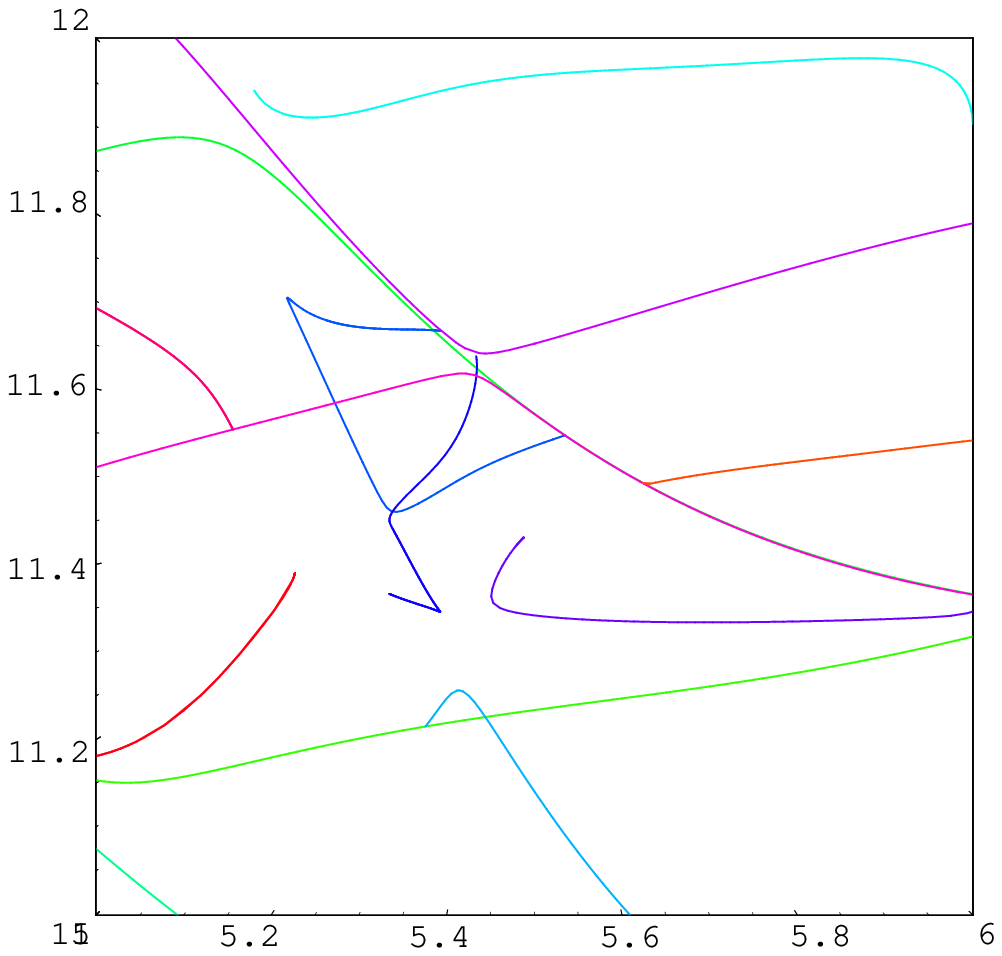}
\hskip .2in \epsffile{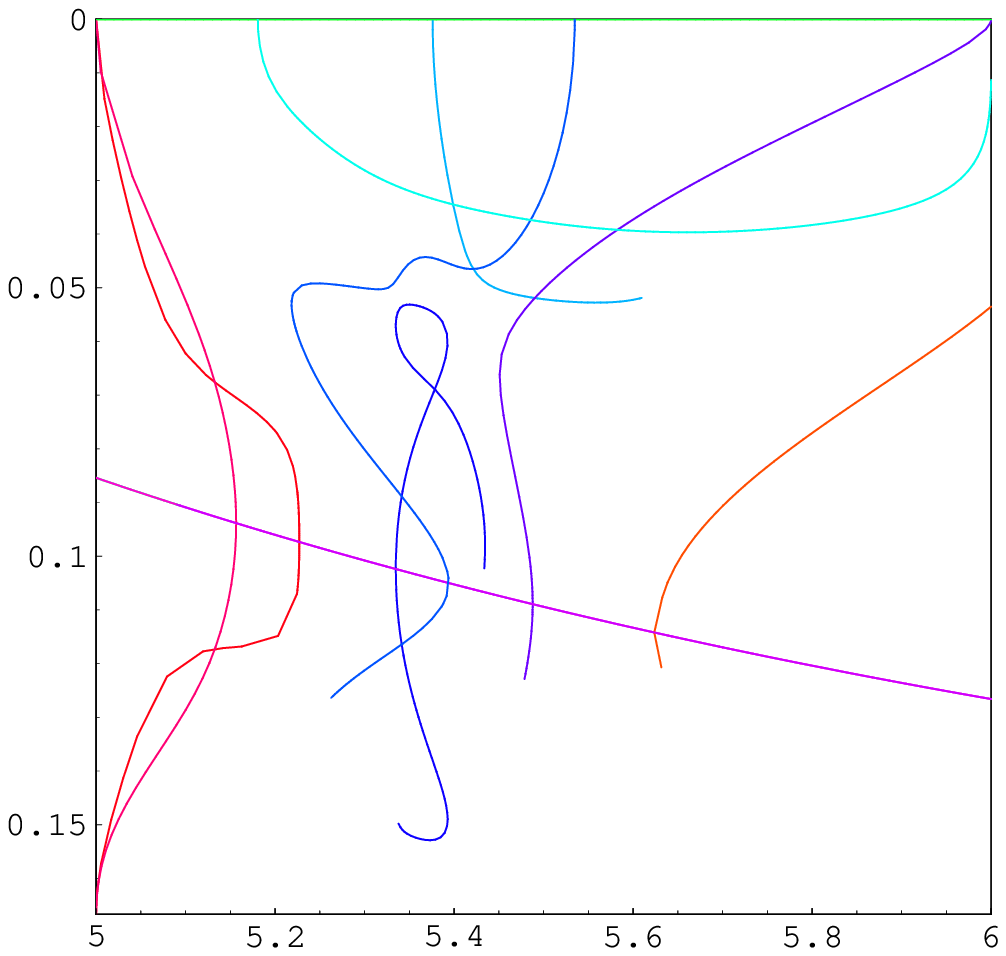}

\centerline{\sl Figure 4.2.1.  Maass forms on deformations of $\Gamma_0^*(5)$.}

\vskip 0.1in

\subhead 4.3 Nonarithmetic examples
\endsubhead

We know that Maass forms
exist on congruence groups,  and 
the deformations of these arithmetic Maass forms give Maass forms on
nonarithmetic groups.  It is natural to ask whether all Maass forms can be
``explained'' by their existence on arithmetic groups.  That is, are all
Maass forms a deformation of a Maass form on an arithmetic group?
We give four examples below which show that, unfortunately, the
answer seems to be ``no.''

\noindent \epsffile{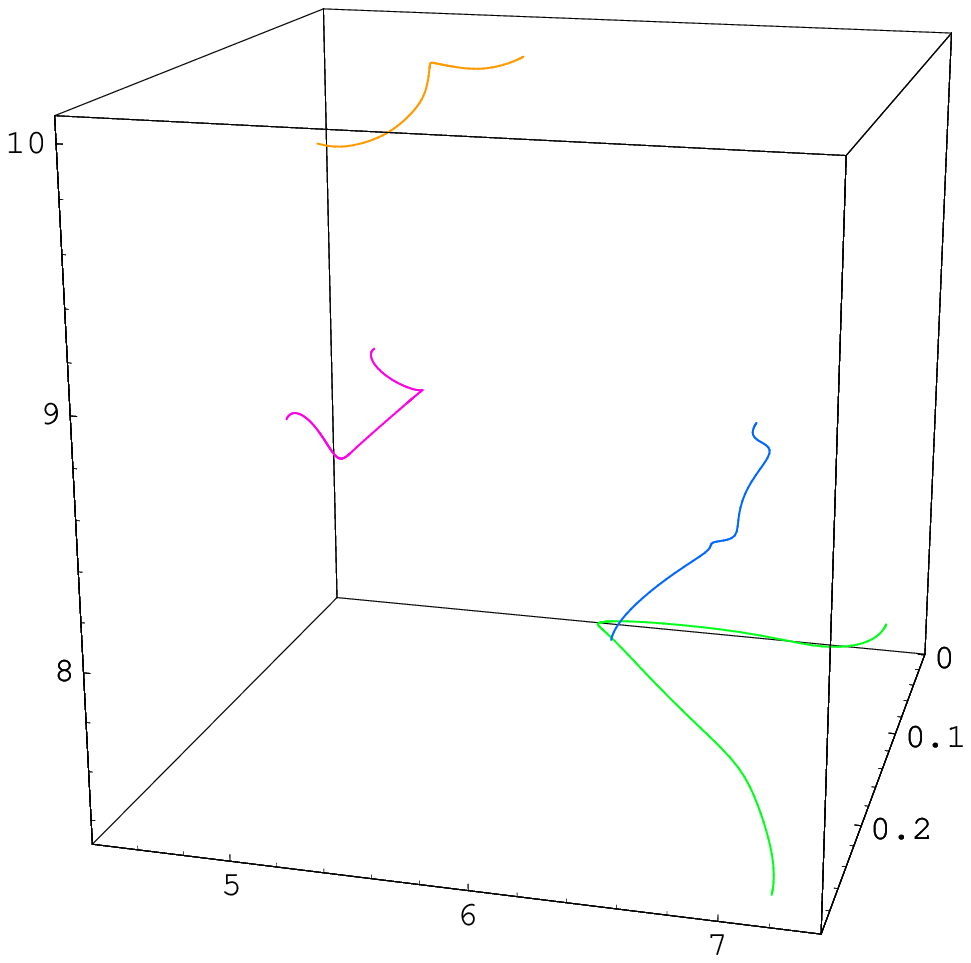}
\hskip .35in \epsffile{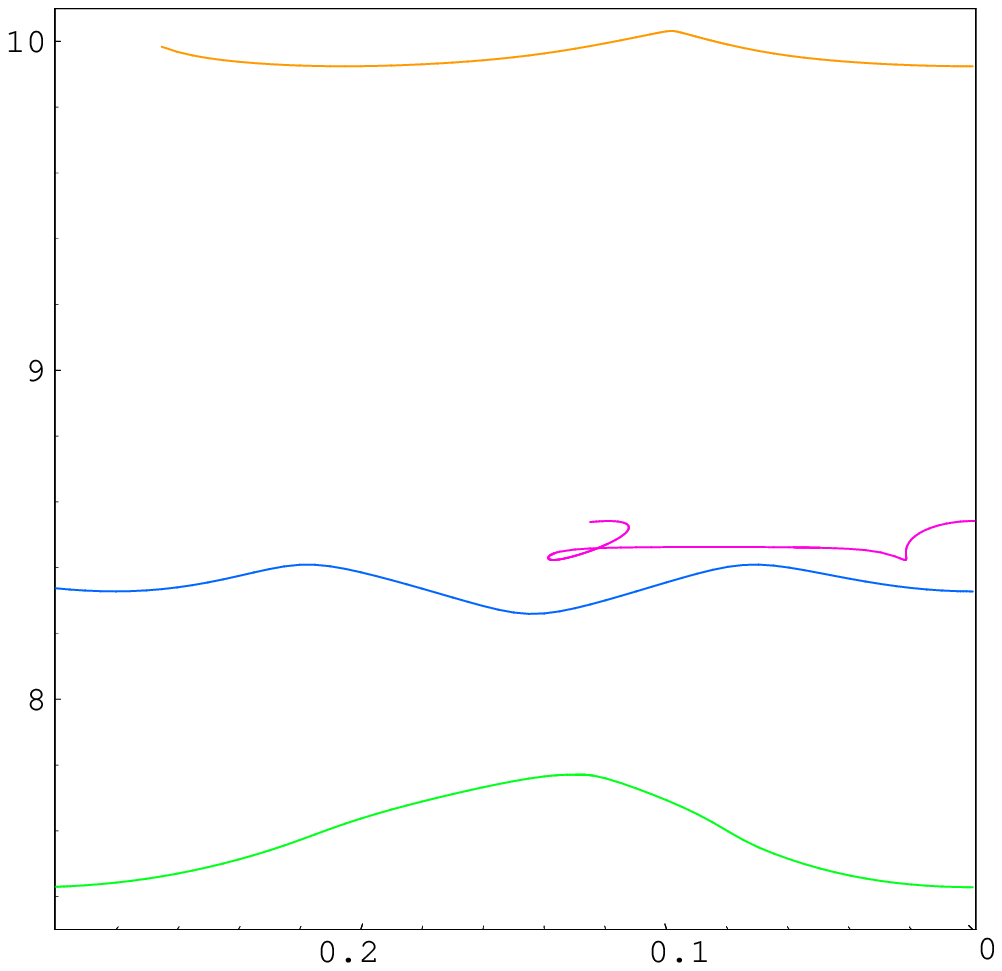}

\noindent \epsffile{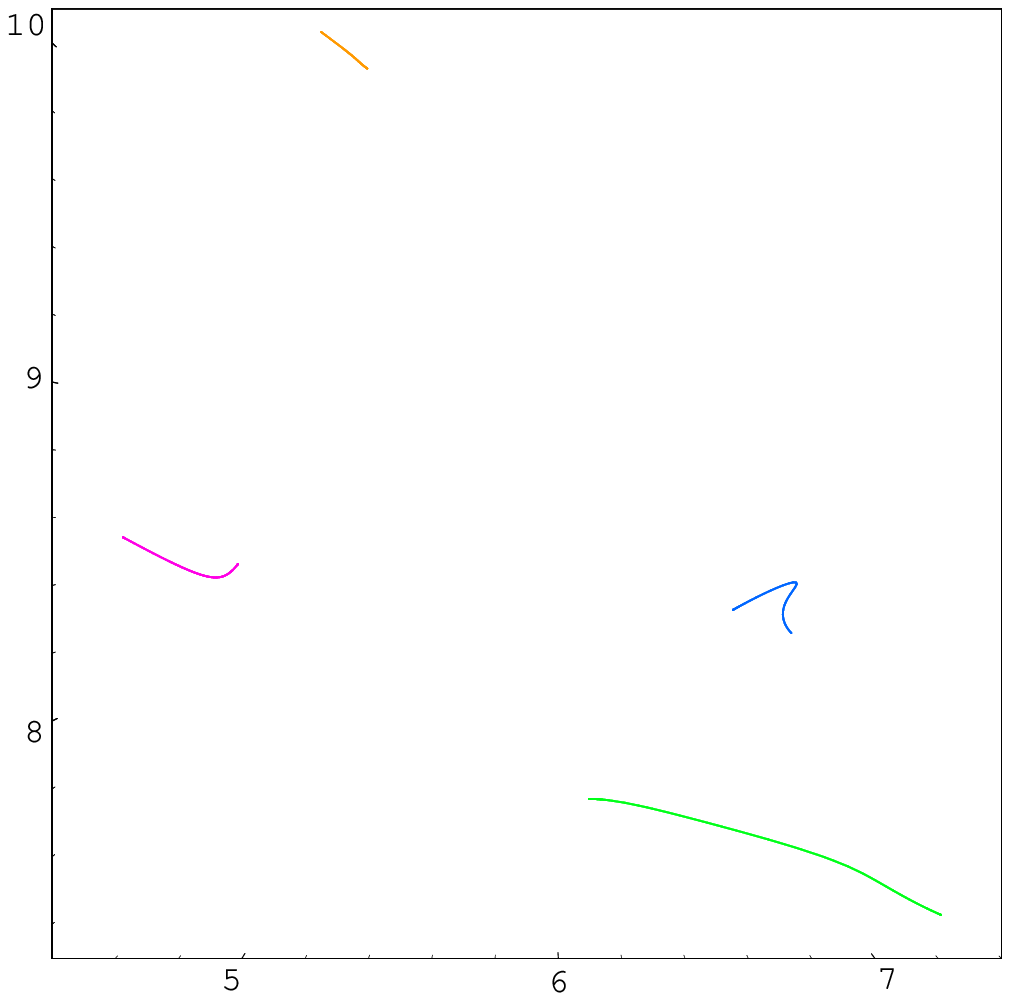}
\hskip .2in \epsffile{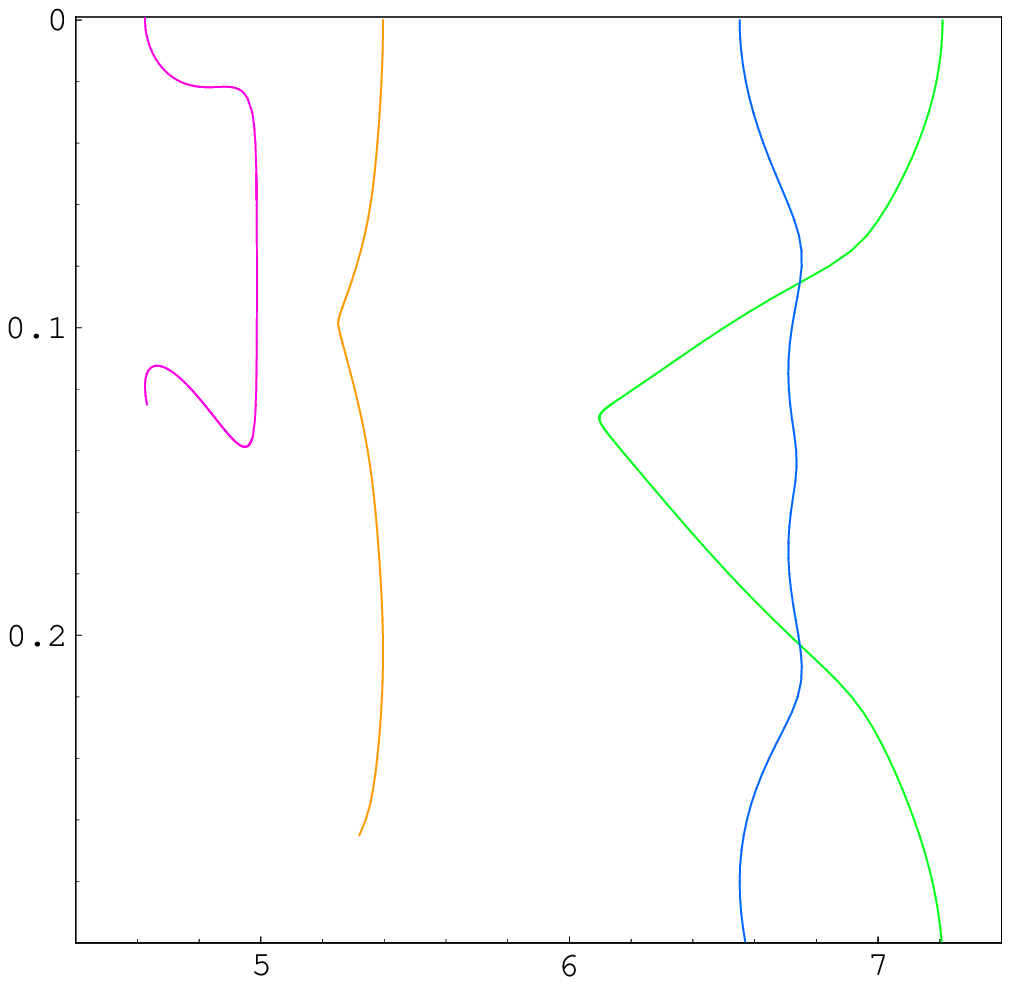}

\centerline{\sl Figure 4.3.1.  Maass forms which do not arise from an arithmetic group}

\vskip 0.1in

Each of the curves in Figure 4.3.1 begins and ends on equivalent groups, so
they correspond to closed paths in Teichm\"uller space.
We believe that none of the points on those curves correspond to 
arithmetic groups.  (We have not been able to verify this directly,
but we have checked that none of those Maass forms have multiplicative 
coefficients.)

\subhead 4.4 Level dynamics of odd forms
\endsubhead

The odd Maass forms on $\Gamma_0(5)$ or $\Gamma_0(6)$ can be
deformed in the $b=0$ plane, which makes their behavior under
deformation somewhat easier to visualize.  Figure 4.4.1 shows
the deformations of the odd Maass forms on $\Gamma_0^*(5)$
for $0<R<16.3$.  The rectangle in the figure corresponds to the
back face of the cube in Figure~4.2.1.   For $a>8$ the deformation 
repeats according to the rule $a\to 4a/(a-4)$.  The limit $a\to 4$
corresponds to the merging of two cusps (and so $a=4$ can never be
reached).  Nevertheless, $\Gamma_{2,2,2}(4,0)$ can be interpreted as
$\Gamma_0^*(4)$, and as
$a\to 4$ the eigenvalues approach those on $\Gamma_0^*(4)$.

Figure~4.2.1 also shows the phenomenon of ``avoided crossing'', 
a manifestation of the fact that a multiplicity in the 
spectrum is a codimension 2 condition.  A closeup of the 
avoided crossing near $(a,R)=(4.83,14.55)$ is given in Figure 4.4.2.
In a future paper we will include a more detailed study of 
the level dynamics of the odd Maass forms.

{

\hskip .05in \epsffile{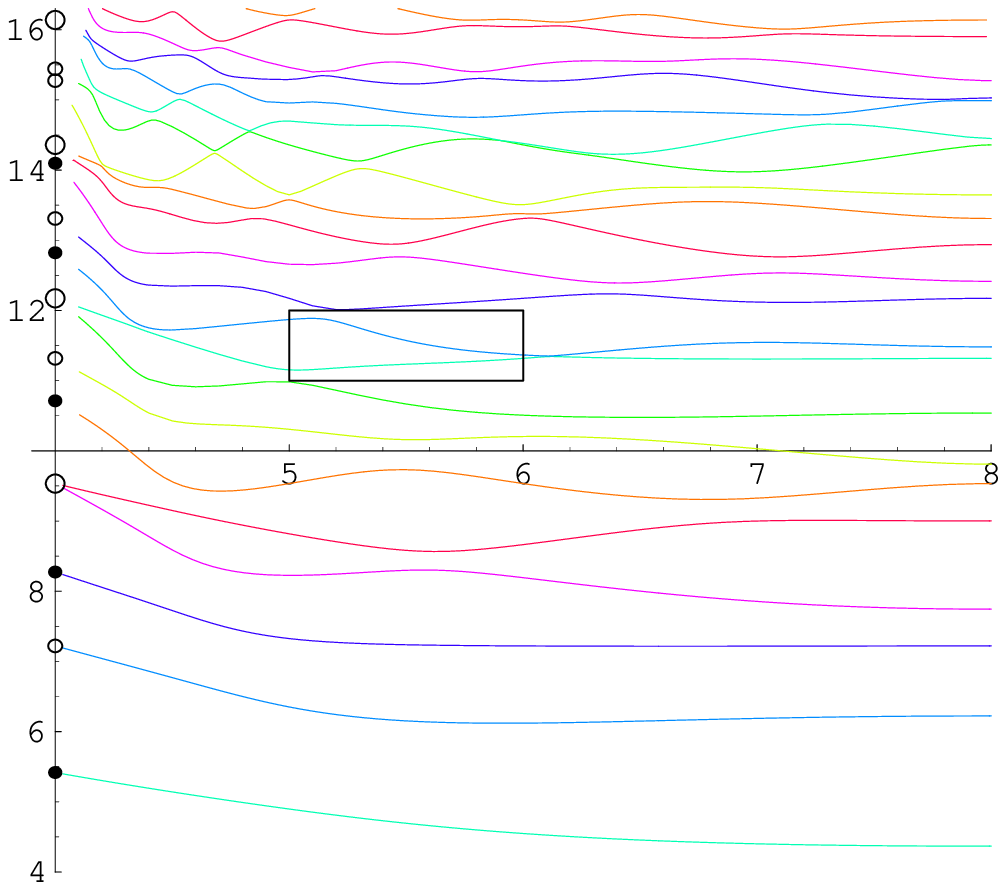}
}

\centerline{\sl Figure 4.4.1.  Odd Maass forms in the $b=0$ plane,
for $4< a\le 8$ and $0\le R\le 16.3$.}
\centerline{\sl The black dots are odd eigenvalues for newforms
on $\Gamma^*_0(4)$, the small circles are} 
\centerline{\sl odd newforms on $\Gamma_0^*(2)$, the large circles are odd eigenvalues
on $\Gamma(1)$, 
 and the}
\centerline{\sl rectangle is the back face of the cube in Figure 4.2.1.}

{

\hskip .35in \epsffile{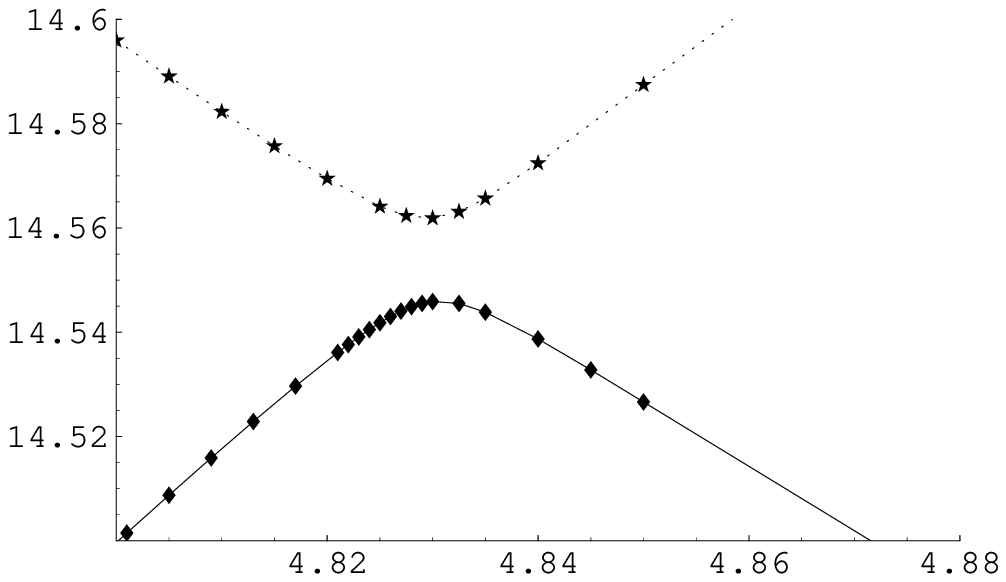}
}

\centerline{\sl Figure 4.4.2.  A closeup of the
avoided crossing near $(a,R)=(4.83,14.55)$.}
\centerline{\sl The dots show the points calculated by the
computer program.}

\subhead 4.5 Deformations of $\Gamma_{2,2,2,2}$
\endsubhead

Since the Teichm\"uller space of $\Gamma_0(11)$ is 4-dimensional, we are not able
to draw pictures like those above.  Our argument that the Maass forms live on
a 3-dimensional set is as follows.

Suppose $f$ is a Maass form on $\Gamma_{2,2,2,2}(a_0,b_0,c_0,d_0)$ with
eigenvalue~$R_0$.  Now suppose we change three of $a_0,b_0,c_0,d_0$ very
slightly.  If the Maass form lives on a 3-dimensional space we will
be able to change $R_0$ and the other parameter slightly and find
another Maass form.  Suppose we make several choices for three of the
parameters (and change them in a slightly different way each time),
and in each case we can slightly adjust $R_0$ and the other 
parameter and again find a Maass form. Then that is persuasive evidence
that the Maass form lives on a 3-dimensional set.

The following table is a representative example of these calculations.
The top line is the data for an odd Maass form of $\Gamma_0(11)$,
which corresponds to $(a,b,c,d)=(-\frac13,0,\frac13,\frac1{2\sqrt{11}})$.
The other entries are deformations of that Maass form, where in each case 
the values given to 5 or fewer decimal places were chosen exactly,
and the values given to 7 decimal places are approximations determined in the search. 

\vskip 0.1in

\hskip 0.5in\vbox{
\+&\phantom{ppppppp} $(a,b,c,d)$ \phantom{pppppppppppppp}& \phantom{ppppppppppp}$R$ \cr
\+&(-0.3333333, 0, 0.3333333, 0.1507556)     & \phantom{ppppppp}  11.8005163 \cr
\+&(-0.31, 0.03, 0.37, 0.1406783) &\phantom{ppppppp}  11.8076532   \cr
\+&(-0.31, 0.03, 0.3704, 0.1404517) & \phantom{ppppppp}  11.8092074   \cr
\+&(-0.31, 0.03, 0.37015, 0.1405936) & \phantom{ppppppp}   11.8082495   \cr
\+&(-0.31, 0.03, 0.3711848, 0.14) &\phantom{ppppppp}  11.8120386   \cr
\+&(-0.3121659, 0.03, 0.37, 0.14) &\phantom{ppppppp}  11.8172635   \cr
}
  
\vskip 0.1in

\subhead 4.6 Additional observations
\endsubhead

Here we address various observations which have not been discussed
previously.

\roster
\item{} It was noted that odd Maass forms on $\Gamma_{2,2,2}(a,b)$ live 
on the plane of deformations with $b=0$.
The figures in the previous sections strongly suggest that even 
Maass forms on $\Gamma_{2,2,2}(a,0)$ live on
a set which meets the $b=0$ plane perpendicularly.
This is confirmed by  Avelin's\cite{A} observation that under deformations
of the form $(a,b)=(5,t)$, the zeros of the scattering
term $\varphi(s)$ have 4th order contact with the $\Re(s)=\frac12$ axis.
The 4th order contact indicates that the deformation is in the (unique)
direction in which the Phillips-Sarnak condition for destruction is not satisfied.
As a specific example, consider the eigenvalue 
$R_0\approx 5.436180461$ for $\Gamma_0(5)$.
Fitting to our data we find that the initial path of the 
deformation is along the curve
$a\approx 5+66.729 b^2-2.1\times 10^4 b^4$ with 
$R\approx R_0 -25.429 b^2+8.3 \times 10^3 b^4$.
Deformation along the path $(a,b)=(5,0) + t (0,1) + t^2 (66.729,0)$
should find that the zero of the scattering matrix has 6th order contact with the
$\Re(s)=\frac12$ axis.  Using the next term 
(assuming our value of $-2.1\times 10^4$
is sufficiently accurate), should show to 8th order contact.

\item{}  There are Maass forms on $\Gamma_{2,2,2}(a,b)$
which can be deformed in two different directions. 
For example, consider the points $(5.53487,0,11.54704)$
and $(5.120,0.0919, 11.9671)$ in Figure 4.2.1.  
The diagrams show two one-parameter deformations
intersecting at each of those points.  
However, there is only one, not two, Maass forms at 
the intersection  point.  That is, those points correspond
to a single Maass form which can be deformed in two independent
directions, as opposed to an eigenvalue of multiplicity two. 
(We ruled out the possibility of a multiple eigenvalue by verifying that 
there were no Maass forms whose first Fourier coefficient vanished.
We also checked that as you approach the point along either path, 
the coefficients have the same limiting value.)

Note that the Maass form at $(5.53487,0,11.54704)$
is odd, so it has the expected deformation in the $b=0$ plane.
However, it also can be deformed along a path which 
initially is perpendicular to the $b=0$ plane.

\item{} There are two places visible in Figure~4.2.1 where the fundamental
domain of the group $\Gamma_{2,2,2}(a,b)$
has an extra symmetry which permits us to distinguish between
``even'' and ``odd'' forms.  One such place is the plane $b=0$, which is the
back face of the cube in the diagrams.  Another place lies along the
curve
from $(a,b)=(5,0.0854)$ to $(a,b)=(6,0.1266)$, which is clearly visible
in Figure~4.2.1.   Since the Eisenstein series are even, the odd Maass forms
have deformations which lie in those regions, which indeed can be seen in
the figure.  

As mentioned in the previous comment, 
the  remaining Maass forms appear to cross those lines of symmetry perpendicularly.
Furthermore, those remaining curves either intersect an arithmetic group,
or they intersect the path of an odd Maass form as they cross one of the lines of symmetry.
So, in some sense, all of these Maass forms have some connection with an arithmetic group.

\item{} The calculations we present here are illuminating, but
the main features are not unexpected. In particular, one should
expect that Maass forms live on $d-1$ dimensional
real analytic subvarieties of Teichm\"uller space.  
As described in Section 1.3, for each cusp
form Phillips and Sarnak 
define a map from Teichm\"uller space to the half-plane~$\Re(s)\le \frac12$,
such that a cusp form is not destroyed if the point lies on the
line~$\Re(s)=\frac12$.  One expects this map
to be real analytic, and
one might also expect it to be non-degenerate, which our
calculations appear to confirm.  Thus, the curves in the above diagrams
are the inverse image of the line $\Re(s)=\frac12$.
It is unclear whether it is surprising that
those curves have
multiple points

Note that the above paragraph concerns the variety on which the 
Phillips and Sarnak integral vanishes, as opposed to the variety
on which the deformed cusp form lives.  However, both varieties
have (real) codimension~1, suggesting that they are the same. 
It is not clear how this relates to Avelin's \cite{A}
observation of scattering zeros having 4th order contact with 
the $\Re(s)=\frac12$ axis.

Sarnak has suggested to us that the dissolving condition for the deformation
of a Maass form can be used to write down a differential equation which
is satisfied by the curves we have found.  It would be interesting to write down
such an equation.

\endroster

\head 5. Questions
\endhead

0. Do the calculations in this paper actually find deformed cusp forms,
or merely Eisenstein series with uniformly very small constant term?
In the latter case, why do the constant terms stay so small?

The remaining questions assume that our calculations actually find
cusp forms.

1. Is there a Weyl's law for equivalence classes of Maass forms?

2.  Can a Maass form ever be deformed to give a different Maass form on 
the same group?  Does the answer change if one restricts
to smooth deformations (so, for example, odd Maass forms in the examples above
can never leave the plane $b=0$)?

3.  In the case of $\Gamma_{2,2,2}(a,b)$, or any other group with
a two real parameter deformation space, how can one detect Maass
forms which can be deformed in two independent directions?
Do those Maass forms have any special properties?

4.  Is it surprising that there exist Maass forms which cannot be smoothly deformed
to give a Maass form on
an arithmetic group?  It is possible that such forms always arise as the
deformation of an ``odd'' Maass form?  

\Refs

\item{[A]} H. Avelin, Research Announcement on the deformation of cusp forms,
U.U.D.M. Report 2002:26, Uppsala Univ.

\item{[Co]} H. Cohn, A numerical survey of the reduction of modular
curve genus by Fricke's involutions, Number theory (New York, 1989/1990), 85-104,
Springer, New York, 1991.

\item{[C1]} Y. Colin de Verdiere, Pseudo-Laplacians I, Ann. Inst. Fourier 
32 (1983), 275-286.

\item{[C2]} Y. Colin de Verdiere, Pseudo-Laplacians II, Ann. Inst. Fourier 
33 (1983), 87-113.

\item{[FJ]} D. Farmer and S. Lemurell, in preparation.

\item{[H1]} D. Hejhal, {\it Eigenvalues of the Laplacian for Hecke triangle groups},
Mem. Amer. Math. Soc. (1992), no. 469, 165pp.

\item{[H2]} D. Hejhal, On Eigenfunctions of the Laplacian for Hecke triangle groups,
in 
{\it Emerging applications of number theory}, Springer, 1999, 291-315.

\item{[H3]} D. Hejhal and S. Arno, On Fourier coefficients of Maass waveforms for
$PSL(2, Z)$,  Math. Comp. 61 (1993), 245-267 and S11-S16.

\item{[Iw]} H. Iwaniec, {\it An introduction to the spectral theory of 
automorphic forms},
Bibl. Rev. Mat. Iber., Madrid, 1995. (Reprinted by AMS).

\item{[L]} W. Luo, Nonvanishing of $L$-values and the Weyl law,
Ann. of Math. (2), 154 (2001) no.~2, 477-502.

\item{[P]} Y. Petridis, Perturbation of scattering poles for hyperbolic
surfaces and central values of $L$-series, Duke Math. J. 103 no 1 (2000), 101-130.

\item{[PS1]} R. Phillips and P. Sarnak, On cusp forms for cofinite
subgroups of $SL(2,R)$, Invent. Math. 80 (1985), 339-364.

\item{[PS2]} R. Phillips and P. Sarnak, The Weyl theorem and the 
deformation of discrete groups, Comm. pure and applied math. 38 (1985), 853-866.

\item{[PS3]} R. Phillips and P. Sarnak, Perturbation theory for the Laplacian
on automorphic functions, J. Amer. Math. Soc. 5 (1992), 1-32.

\item{[S1]}  P. Sarnak, On cusp forms.  {\it The Selberg trace formula and 
related topics (Brunswick, Maine, 1984)},  393--407, 
Contemp. Math., 53, Amer. Math. Soc., Providence, RI, 1986

\item{[S2]}  P. Sarnak, On cusp forms. II.  {\it Festschrift in honor of 
I. I. Piatetski-Shapiro on the occasion of his sixtieth birthday, Part II 
(Ramat Aviv, 1989)},  237--250, 
Israel Math. Conf. Proc., 3, Weizmann, Jerusalem, 1990

\item{[Ta]} K. Takeuchi, A characterization of arithmetic Fuchsian
groups, J. Math. Soc. Japan, 27 (1975), 600-612.

\item{[W]} S. Wolpert, Disappearance of cusp forms in special families,
Ann. of Math. (2) 139 (1994) no.~2, 239-291.

\endRefs

\enddocument

\bye